\documentclass{article}
\usepackage{amsfonts}
\usepackage{amsmath}
\usepackage{amsmath,amssymb,latexsym,color}
\usepackage[mathscr]{eucal}

\newtheorem{theorem}{Theorem}
\newtheorem{prop}{Proposition}
\newtheorem{lemma}{Lemma}
\newtheorem{rem}{Remark}
\newtheorem{cor}{Corollary}

\newtheorem{fact}{Fact}

\newtheorem{defin}{Definition}

\newtheorem{remark}{Remark}[section]

\newcommand{\proof}{{\noindent\it Proof.\quad}}
\newcommand{\qed}{\hfill\Box\medskip}
\begin{document}

\title
{\bf Transformations of polar Grassmannians preserving certain
intersecting relations}

\author{
Wen Liu$^{\rm a,b}$\quad Mark Pankov$^{\rm c}$ \quad Kaishun
Wang$^{\rm a}$\footnote{Corresponding author.\newline  Email
addresses: liuwen1975@126.com (Wen Liu), pankov@matman.uwm.edu.pl
(Mark
Pankov), wangks@bnu.edu.cn (Kaishun Wang)}  \\
{\footnotesize  $^{\rm a}$ \em  Sch. Math. Sci. {\rm \&} Lab. Math.
Com. Sys.,
Beijing Normal University, Beijing, 100875,  China}\\
\footnotesize $^{\rm b}$ \em  Ins.  Math. {\rm \&}  Inf., Hebei
Normal University, Shijiazhuang, 050024, China\\
\footnotesize $^{\rm c}$ \em Department of Mathematics and Computer
Science, University of Warmia and Mazury Olsztyn, Poland}

\date{}
\maketitle
\begin{abstract}
Let $\Pi$ be a polar space of rank $n\ge 3$. Denote by ${\mathcal
G}_{k}(\Pi)$ the polar Grassmannian formed by singular subspaces of
$\Pi$ whose projective dimension is equal to $k$. Suppose that $k$
is an integer not greater than $n-2$ and consider the relation
${\mathfrak R}_{i,j}$, $0\le i\le j\le k+1$ formed by all pairs
$(X,Y)\in {\mathcal G}_{k}(\Pi)\times {\mathcal G}_{k}(\Pi)$ such
that $\dim_{p}(X^{\perp}\cap Y)=k-i$ and $\dim_{p} (X\cap Y)=k-j$
($X^{\perp}$ consists of all points of $\Pi$ collinear to every
point of $X$). We show  that every bijective transformation of
${\mathcal G}_{k}(\Pi)$ preserving ${\mathfrak R}_{1,1}$ is induced
by an automorphism of $\Pi$ and the same holds for the relation
${\mathfrak R}_{0,t}$ if $n\ge 2t\ge 4$ and $k=n-t-1$. In the case
when $\Pi$ is a finite classical polar space, we establish that the
valencies of ${\mathfrak R}_{i,j}$ and ${\mathfrak R}_{i',j'}$ are
distinct if $(i,j)\ne (i',j')$.



\end{abstract}

\section{Introduction}
Let $V$ be an $n$-dimensional vector space $V$ (over a division ring).
Denote by ${\mathcal G}_{k}(V)$ the Grassmannian formed by $k$-dimensional subspaces of $V$.
Suppose that $1<k\le n-k$.
For every integer $i$ satisfying $1\le i\le k$ we define the following relation
$${\mathfrak R}_{i}:=\{(X,Y)\in{\mathcal G}_{k}(V)\times{\mathcal G}_{k}(V): \dim_{l}(X\cap Y)=k-i\}.$$
We write $\dim_{l}$ for the linear dimension (the dimension of vector spaces and their subspaces),
since we want to distinguish it from the projective dimension
(the dimension of projective spaces and their subspaces) which will be denoted by $\dim_{p}$.
Note that $(X,Y)\in {\mathfrak R}_{i}$ if and only if
the distance between $X$ and $Y$ in the {\it Grassmann graph}
$$\Gamma_{k}(V)=({\mathcal G}_{k}(V), {\mathfrak R}_{1})$$
is equal to $i$.

The same relations can be defined on dual polar spaces.
Let $\Pi$ be a polar space of rank $n$.
Denote by ${\mathcal G}_{k}(\Pi)$ the polar Grassmannian formed by
singular subspaces of $\Pi$ whose projective dimension is equal to $k$.
The associated dual polar space is formed by maximal singular subspaces,
i.e. singular subspaces of dimension $n-1$.
For any integer $i$ satisfying $1\le i\le n$ we define
$${\mathfrak R}_{i}:=\{(X,Y)\in{\mathcal G}_{n-1}(\Pi)\times{\mathcal G}_{n-1}(\Pi): \dim_{p}(X\cap Y)=n-1-i\}.$$
As above, we have $(X,Y)\in {\mathfrak R}_{i}$ if and only if
the distance between $X$ and $Y$ in the {\it dual polar graph}
$$\Gamma_{n-1}(\Pi)=({\mathcal G}_{n-1}(\Pi), {\mathfrak R}_{1})$$
is equal to $i$.

By \cite{Chow}, every automorphism of the Grassmann graph $\Gamma_{k}(V)$
is induced by a semilinear automorphism of $V$ or a semilinear isomorphism of $V$ to
the dual vector space $V^{*}$ (the second possibility can be realized only in the case when $n=2k$).
Similarly, every automorphism of the dual polar graph $\Gamma_{n-1}(\Pi)$
is induced by an automorphism of the polar space $\Pi$.
The latter was proved by Chow \cite{Chow} for classical polar spaces only,
but Chow's method works in the general case \cite[Section 4.6]{Pankov-book}.
Some results closely related to these statements were obtained
\cite{Havlicek1,Havlicek2,HavlicekPankov,Huang1,Huang2,Huang3,Kreuzer,Pankov1,Pankov2}
and we refer \cite{Pankov-book} for a survey.

Every bijective transformation of ${\mathcal G}_{k}(V)$ preserving ${\mathfrak R}_{k}$
is an automorphism of $\Gamma_{k}(V)$ \cite{BH,HavlicekPankov}.
For the relation ${\mathfrak R}_{i}$ with $1<i<k$ the same is not proved.
However, all bijective transformations of ${\mathcal G}_{k}(V)$ preserving
${\mathfrak R}_{1}\cup\dots\cup {\mathfrak R}_{m}$
are automorphisms of $\Gamma_{k}(V)$ for every integer $m<k$ \cite{Lim}.
This is a generalization of the previous result; indeed, if $m=k-1$ then
the considered above transformations preserve ${\mathfrak R}_{k}$.
The same statement is proved for some dual polar spaces \cite{Huang5}.
Results of similar nature were established for other objects \cite{AbHVM, HVM1, HVM2, HH, Huang4}.

Now suppose that $k$ is an integer not greater than $n-2$
and consider the relation ${\mathfrak R}_{i,j}$, $0\le i\le j\le k+1$
formed by all pairs
$$(X,Y)\in {\mathcal G}_{k}(\Pi)\times {\mathcal G}_{k}(\Pi)$$
satisfying the following conditions
$$\dim_{p}(X^{\perp}\cap Y)=k-i\;\mbox{ and }\;\dim_{p} (X\cap Y)=k-j$$
($X^{\perp}$ consists of all points of $\Pi$ collinear to every point of $X$).
All automorphisms of the {\it polar Grassmann graph}
$$\Gamma_{k}(\Pi)=({\mathcal G}_{k}(\Pi), {\mathfrak R}_{0,1})$$
are described in \cite[Section 4.6]{Pankov-book}.
They are induced by automorphisms of $\Pi$,  except the case when $n=4$, $k=1$
and our polar space is of type $\textsf{D}$
(every $(n-2)$-dimensional singular subspace is contained in precisely two maximal singular subspaces).
Also, every automorphism of so-called {\it weak Grassmann graph}
$$\Gamma^{w}_{k}(\Pi)=({\mathcal G}_{k}(\Pi),{\mathfrak R}_{0,1}\cup {\mathfrak R}_{1,1})$$
is induced by an automorphism of $\Pi$ \cite{PPZ}.

We show  that every bijective transformation of ${\mathcal
G}_{k}(\Pi)$ preserving ${\mathfrak R}_{1,1}$ is induced by an
automorphism of $\Pi$ (Theorem \ref{theorem1}). Our second result
(Theorem \ref{theorem2}) states that the same holds for the relation
${\mathfrak R}_{0,t}$ if $n\ge 2t\ge 4$ and $k=n-t-1$. The latter
relation is equivalent to the fact that two elements with the
relation ${\mathfrak R}_{0,t}$ span a maximal singular subspace.

Note that for   finite symplectic and hermitian polar spaces,
the first result under some conditions was proved in
\cite{liu,zeng}.

For a finite classical polar space, the valency of every
${\mathfrak R}_{i,j}$ was given by Stanton \cite{Stanton}.  We
establish that the valencies of ${\mathfrak R}_{i,j}$ and
${\mathfrak R}_{i',j'}$ are distinct if $(i,j)\ne (i',j')$
(Proposition \ref{prop-val}).

\section{Polar spaces}
We recall  some basic properties of polar spaces
and refer \cite{BC, Pankov-book,Ueberberg} for their proofs.

Let ${\mathcal P}$ be a non-empty set whose elements are called {\it points}
and ${\mathcal L}$ be a family formed by proper subsets of ${\mathcal P}$ called {\it lines}.
Two distinct points joined by a line are said to be {\it collinear}.
Let $\Pi=(P,{\mathcal L})$ be a {\it partial linear space}, i.e.
each line contains at least two points and
for any distinct collinear points $p,q\in {\mathcal P}$ there is precisely one line containing them,
this line is denoted by $pq$.

We say that $S\subset {\mathcal P}$ is a {\it subspace} of $\Pi$
if for any distinct collinear points $p,q\in S$
the line $pq$ is contained in $S$.
A {\it singular} subspace is a subspace where any two distinct points are collinear.
Note that the empty set and a single point are singular subspaces.
Using Zorn lemma, we show that every singular subspace is contained
in a certain maximal singular subspace.

From this moment we suppose that $\Pi$ is a {\it polar space}.
This means that the following axioms hold:
\begin{enumerate}
\item[(P1)] each line contains at least three points,
\item[(P2)] there is no point collinear to all points,
\item[(P3)] if $p\in {\mathcal P}$ and $L\in {\mathcal L}$
then $p$ is collinear to precisely one point or all points of the line $L$,
\item[(P4)] every flag formed by singular subspaces is finite.
\end{enumerate}
If there is a maximal singular subspace of $\Pi$ containing more than one line
then all maximal singular subspaces of $\Pi$ are projective spaces
of the same finite dimension $\ge 2$.
We say that the {\it rank} of $\Pi$ is $n$ if this dimension is equal to $n-1$.

In the case when the rank of $\Pi$ is not less than $4$,
every maximal singular subspace $M$ can be identified with
the projective space associated to a certain $n$-dimensional vector space $V$ (over a division ring).
Then every non-empty singular subspace $S\subset M$ will be identified with
the corresponding subspace of the vector space $V$.

The collinearity relation on $\Pi$ is denoted by $\perp$.
We write $p\perp q$ if $p,q\in {\mathcal P}$ are collinear points and $p\not\perp q$ otherwise.
If $X,Y\subset {\mathcal P}$ then $X\perp Y$ means that every point of $X$ is collinear to all points of $Y$.
For every subset $X\subset {\mathcal P}$ satisfying $X\perp X$
the minimal singular subspace containing $X$ is called {\it spanned} by $X$
and denoted by $\langle X \rangle$.
For every subset $X\subset {\mathcal P}$ we denote by $X^{\perp}$
the subspace of $\Pi$ formed by all points  collinear to all points of $X$.

\begin{fact}\label{fact1}
Let $X$ be a subset of ${\mathcal P}$ satisfying $X\perp X$ and spanning a maximal singular subspace $M$.
Then $p\perp X$ implies that $p\in M$.
\end{fact}

\begin{fact}\label{fact2}
If $M$ is a maximal singular subspace of $\Pi$
then for every point $p\in {\mathcal P}$ such that $p\not\in M$ we have
$$\dim_{p} (p^{\perp}\cap M)=n-2.$$
\end{fact}

\begin{fact}\label{fact3}
For every singular subspace $S$ there are maximal singular subspaces $M$ and $N$
such that  $S=M\cap N$.
\end{fact}

\section{Results}
Let $\Pi=({\mathcal P},{\mathcal L})$ be a polar space of   rank
$n\ge 3$. Recall that the polar Grassmannian formed by all
$k$-dimensional singular subspaces of $\Pi$ is denoted by ${\mathcal
G}_{k}(\Pi)$. Suppose that $k\le n-2$. We consider the relation
${\mathfrak R}_{i,j}$, $0\le i\le j\le k+1$ is formed by all pairs
$$(X,Y)\in {\mathcal G}_{k}(\Pi)\times{\mathcal G}_{k}(\Pi)$$
satisfying
$$\dim_{p} (X^{\perp}\cap Y)=k-i\;\mbox{ and }\; \dim_{p} (X\cap Y)=k-j.$$
Since the equality
$$\dim_{p} (X^{\perp}\cap Y)=\dim_{p} (Y^{\perp}\cap X)$$
holds for any pair $X,Y\in {\mathcal G}_{k}(\Pi)$,
this relation is symmetric.
Every automorphism of $\Pi$ (a bijective transformation of ${\mathcal P}$ preserving ${\mathcal L}$)
induces a transformation of ${\mathcal G}_{k}(\Pi)$ which preserves all ${\mathfrak R}_{i,j}$.

First, we determine all automorphisms of the graph
$$\Gamma'_{k}(\Pi)=({\mathcal G}_{k}(\Pi),{\mathfrak R}_{1,1}).$$

\begin{theorem}\label{theorem1}
Every automorphism of $\Gamma'_{k}(\Pi)$ is induced by an automorphism of $\Pi$.
\end{theorem}

\begin{rem}{\rm
For $k=0$ this statement is trivial.
The edges of $\Gamma'_{0}(\Pi)$ are pairs of non-collinear points of $\Pi$
and every automorphism of this graph is an automorphism of
the collinearity graph $\Gamma_{0}(\Pi)$.
It is well-known that the class of automorphisms of $\Gamma_{0}(\Pi)$ coincides with
the class of automorphisms of $\Pi$.
}\end{rem}

In the case when $k\in \{1,\dots, n-3\}$,
the distance between $S,U\in {\mathcal G}_{k}(\Pi)$ in the Grassmann graph $\Gamma_{k}(\Pi)$
is equal to $2$ if and only if $(S,U)$ belongs to ${\mathfrak R}_{1,1}\cup {\mathfrak R}_{0,2}$.
The distance between $S,U\in {\mathcal G}_{n-2}(\Pi)$ in $\Gamma_{n-2}(\Pi)$
is equal to $2$ if and only if $(S,U)\in{\mathfrak R}_{1,1}$
(if $k=n-2$ then ${\mathfrak R}_{0,2}$ is empty).
Theorem \ref{theorem1} gives the following.

\begin{cor}
Let $f$ be a bijective transformation of ${\mathcal G}_{n-2}(\Pi)$
satisfying the following condition:
the distance between $S,U\in {\mathcal G}_{n-2}(\Pi)$ in $\Gamma_{n-2}(\Pi)$ is equal to $2$
if and only if the distance between $f(S)$ and $f(U)$ in $\Gamma_{n-2}(\Pi)$ is equal to $2$.
Then $f$ is an automorphism of $\Gamma_{n-2}(\Pi)$.
\end{cor}

Suppose that $k=n-t-1$, where $t$ is an integer satisfying $n\ge
2t\ge 4$. Then $(S,U)\in {\mathfrak R}_{0,t}$ is equivalent to the
fact that $S$ and $U$ span a maximal singular subspace of $\Pi$. Our
second result describes all automorphisms of the graph
$$\Gamma''_{k}(\Pi)=({\mathcal G}_{k}(\Pi),{\mathfrak R}_{0,t}).$$

\begin{theorem}\label{theorem2}
Every   automorphism of $\Gamma''_{n-t-1}(\Pi)$ is induced by an
automorphism of $\Pi$.
\end{theorem}

\section{Cliques}
From this moment we suppose that $k\in\{1,\dots,n-2\}$.
For every singular subspace $N$ such that $\dim_{p} N<k$
we denote by $[N\rangle_{k}$ the set of all elements of ${\mathcal G}_{k}(\Pi)$ containing $N$.
This subset is said to be a {\it big star} if $N\in {\mathcal G}_{k-1}(\Pi)$.

If $N$ and $M$ are singular subspaces satisfying
$$N\subset M\;\mbox{ and }\;\dim_{p} N<k<\dim_{p} M$$
then we denote by $[N,M]_{k}$ the set of all $S\in {\mathcal G}_{k}(\Pi)$
such that $N\subset S\subset M$.
This subset is called a {\it star} if
$$N\in {\mathcal G}_{k-1}(\Pi)\;\mbox{ and }\;M\in {\mathcal G}_{n-1}(\Pi).$$
In the case when $N=\emptyset$, we write $\langle M]_{k}$ instead of $[N,M]_{k}$.
We say that $\langle M]_{k}$ is a {\it top} if $M\in {\mathcal G}_{k+1}(\Pi)$.

All maximal cliques of the Grassmann graph
$$\Gamma_{k}(\Pi)=({\mathcal G}_{k}(\Pi),{\mathfrak R}_{0,1})$$
and the weak Grassmann graph
$$\Gamma^{w}_{k}(\Pi)=({\mathcal G}_{k}(\Pi),{\mathfrak R}_{0,1}\cup {\mathfrak R}_{1,1})$$
are known \cite[Subsection 4.5.1 and Subsection 4.6.2]{Pankov-book}.
Every maximal clique of $\Gamma^{w}_{k}(\Pi)$ is a big star or a top.
Every maximal clique of $\Gamma_{k}(\Pi)$ is a star or a top.
In the case when $k=n-2$, every star is contained in a certain top
and all maximal cliques of $\Gamma_{k}(\Pi)$ are tops.

\begin{prop}\label{prop-cliq1}
Every clique of  $\Gamma'_{k}(\Pi)$ is contained in a big star.
\end{prop}

\proof If ${\mathcal C}$ is a clique of $\Gamma'_{k}(\Pi)$ then it
is a clique of $\Gamma^{w}_{k}(\Pi)$. Hence ${\mathcal C}$ is
contained in a big star or a top. Since any two distinct elements of
a top are non-adjacent vertices of $\Gamma'_{k}(\Pi)$, ${\mathcal
C}$ is a subset in a big star. $\qed$

Let $N\in {\mathcal G}_{k-1}(\Pi)$.
For every $M\in {\mathcal G}_{k+1}(\Pi)$ containing $N$
the subset $[N,M]_{k}$ is said to be a {\it line}.
The big star $[N\rangle_{k}$ together with all lines defined above is a polar space of rank $n-k$
\cite[Lemma 4.4]{Pankov-book}.
This polar space will be denoted by $\Pi_{N}$.

\begin{lemma}\label{lemma-cliq1}
Let $N\in {\mathcal G}_{k-1}(\Pi)$.
Two distinct elements of the big star $[N\rangle_{k}$
are adjacent vertices of the graph $\Gamma'_{k}(\Pi)$ if and only if
they are non-collinear points of the polar space $\Pi_{N}$.
\end{lemma}

\proof Easy verification. $\qed$

\begin{prop}\label{prop-cliq2}
Suppose that $k=n-t-1$, where $t$ is an integer satisfying $n\ge
2t\ge 4$. If $S$ and $U$ are adjacent vertices of
$\Gamma''_{k}(\Pi)$, then $M:=\langle S\cup U\rangle$ is a maximal
singular subspace of $\Pi$ and every clique of $\Gamma''_{k}(\Pi)$
containing $S$ and $U$ is a subset in $\langle M]_{k}$.
\end{prop}

\proof It is clear that $S\perp U$ and $M$  is a maximal singular
subspace. Let ${\mathcal C}$ be a clique of $\Gamma''_{k}(\Pi)$
containing $S$ and $U$. Then for every $A\in {\mathcal C}$ we have
$A\perp S$ and $A\perp U$. By Fact \ref{fact1}, this implies that
$A\subset M$. Hence ${\mathcal C}$ is contained in $\langle M]_{k}$.
$\qed$

\section{Proof of Theorem \ref{theorem1}}

\begin{lemma}\label{lemma1-1}
Let $p$ and $q$ be non-collinear points of $\Pi$
and let $t$ be a point of $\Pi$ collinear to at least one of the points  $p,q$.
Then there exists a point of $\Pi$ non-collinear to $p,q,t$.
\end{lemma}

\proof Suppose that  the statement fails and every point of $\Pi$ is
collinear to at least one of the points $p,q,t$, in other words,
\begin{equation}\label{eq1}
{\mathcal P}=p^{\perp}\cup q^{\perp}\cup t^{\perp}.
\end{equation}
First we show that every maximal singular subspace of $\Pi$
contains at least one of the points $p,q,t$.

Let $M$ be a maximal singular subspace of $\Pi$.
If each of the points $p,q,t$ does not belong to $M$ then
$$p^{\perp}\cap M,\;q^{\perp}\cap M,\; t^{\perp}\cap M$$
are $(n-2)$-dimensional subspaces of $M$ (Fact \ref{fact2})
and \eqref{eq1} implies that $M$ is the union of these subspaces.
The latter is impossible,
since a projective space cannot be presented as the union of three hyperplanes.

By our assumption, $t$ is collinear to at least one of the points  $p$ and $q$.
Suppose that $q\perp t$.
Since $p\not\perp q$,
the line $qt$ contains the unique point $s$ collinear to $p$.
One of the following possibilities is realized:
\begin{enumerate}
\item[(1)] $s=t$,
\item[(2)] $s\ne t$.
\end{enumerate}
In the case (1), we take any point $v$ on the line $qt$ different from $q$ and $t$.
It is clear that $p\not\perp v$.
By Fact \ref{fact3}, there exist maximal singular subspaces $M$ and $N$ such that
$M\cap N=\{v\}$.
Since $p\not\perp v$, they do not contain $p$.
Then one of these subspaces contains $q$ and the other contains $t$.
So, each of these subspaces contains two distinct points of the line $qt$.
This means that this line is contained in $M\cap N$ which is impossible.

In the case (2), we take any point $w$ on the line $ps$ different
from $p$ and $s$. As in the previous case, we consider maximal
singular subspaces $M$ and $N$ such that $M\cap N=\{w\}$. One of
these subspaces contains $q$ or $t$. Then at least one of the points
$q,t$ is collinear to $w$. Since $q$ and $t$ both are collinear to
$s$ and $w$ is on the line $ps$, one of the points $q,t$ is
collinear to all points of the line $ps$. Thus $p$ is collinear to
$q$ or $t$  which is impossible. $\qed$

\begin{lemma}\label{lemma1-2}
Suppose that $N\in {\mathcal G}_{k-1}(\Pi)$.
Let $P,Q\in [N\rangle_{k}$ be adjacent vertices of $\Gamma'_{k}(\Pi)$
and let $T\in [N\rangle_{k}$ be a vertex of $\Gamma'_{k}(\Pi)$
non-adjacent to at least one of the vertices $P,Q$.
Then there exists $S\in [N\rangle_{k}$ adjacent to $P,Q,T$ in $\Gamma'_{k}(\Pi)$.
\end{lemma}

\proof By Lemma \ref{lemma-cliq1}, $P$ and $Q$ are non-collinear
points of $\Pi_{N}$ and $T$ is a point of $\Pi_{N}$ collinear to at
least one of the points $P,Q$. We apply Lemma \ref{lemma1-1} to the
polar space $\Pi_{N}$ and get the claim. $\qed$

Let $f$ be an automorphism of $\Gamma'_{k}(\Pi)$.

Show that $f$ transfers big stars to subsets of big stars.
We take any $N\in {\mathcal G}_{k-1}(\Pi)$.
Let $P$ and $Q$ be adjacent vertices of $\Gamma'_{k}(\Pi)$
contained in the big star $[N\rangle_{k}$.
Then $f(P)$ and $f(Q)$ are adjacent vertices of $\Gamma'_{k}(\Pi)$
contained in the big star $[N'\rangle_{k}$, where
$$N'=f(P)\cap f(Q).$$
We assert that $f(T)\in [N'\rangle_{k}$ for every $T\in [N\rangle_{k}$
and prove this statement in several steps.

(i). First, we consider the case when
$T$ is a vertex of  $\Gamma'_{k}(\Pi)$ adjacent to both $P$ and $Q$.
Then $f(P),f(Q),f(T)$
form a clique in $\Gamma'_{k}(\Pi)$ which,
by Proposition \ref{prop-cliq1}, is contained in a certain big star $[N''\rangle_{k}$.
We have
$$N'=f(P)\cap f(Q)=N''$$
which gives the claim.

(ii). Consider the case when $T$ is a vertex of  $\Gamma'_{k}(\Pi)$
adjacent to precisely one of the vertices $P,Q$.
Suppose that $T$ is adjacent to $P$.
Lemma \ref{lemma1-2} implies the existence of
a vertex $S\in [N\rangle_{k}$ in the graph $\Gamma'_{k}(\Pi)$ adjacent to $P,Q,T$.
By (i), $f(S)$ belongs to $[N'\rangle_{k}$.
Then
$$f(P)\cap f(S)=N'$$
and $f(P),f(S),f(T)$ form a clique of $\Gamma'_{k}(\Pi)$.
As in (i), we show that $f(T)\in [N'\rangle_{k}$.

(iii). Suppose that $T$ is a vertex of $\Gamma'_{k}(\Pi)$ non-adjacent to both $P$ and $Q$.
As above, we consider $S\in [N\rangle_{k}$ which is a vertex of $\Gamma'_{k}(\Pi)$ adjacent to $P,Q,T$
and obtain that $f(S)\in [N'\rangle_{k}$.
We apply the arguments from (ii) to $P,S,T$ and establish that $f(T)\in [N'\rangle_{k}$.

So, $f$ transfers big stars to subsets of big stars.
The same arguments show that $f^{-1}$ sends big stars to subsets of big stars.
This means that $f$ and $f^{-1}$ both map big stars to big stars, i.e.
there exists a bijective transformation $g$ of ${\mathcal G}_{k-1}(\Pi)$
such that
$$f([N\rangle_{k})=[g(N)\rangle_{k}$$
for every $N\in {\mathcal G}_{k-1}(\Pi)$.

Let $U\in {\mathcal G}_{k}(\Pi)$.
Then $g$ transfers the top $\langle U]_{k-1}$ to the top $\langle f(U)]_{k-1}$.
Indeed, we have
$$N\in \langle U]_{k-1}\Leftrightarrow U\in [N\rangle_{k}\Leftrightarrow f(U)\in [g(N)\rangle_{k}\Leftrightarrow
g(N)\in \langle f(U)]_{k-1}.$$
Similarly, $g^{-1}$ sends $\langle U]_{k-1}$
to the top $\langle f^{-1}(U)]_{k-1}$.

Therefore, $g$ and $g^{-1}$ both transfer tops to tops.
Since for any two adjacent vertices of the Grassmann graph $\Gamma_{k-1}(\Pi)$ there is a top containing them,
$g$ is an automorphism of $\Gamma_{k-1}(\Pi)$.
In some special cases,
there are automorphisms of $\Gamma_{k-1}(\Pi)$ which are not induced by automorphisms of $\Pi$.
However, every automorphism of $\Gamma_{k-1}(\Pi)$ transferring tops to tops is
induced by an automorphism of $\Pi$ \cite[Section 4.6.1]{Pankov-book}.
Thus $g$ is induced by an automorphism of $\Pi$.
An easy verification shows that this automorphism also induces $f$.

\section{Proof of Theorem \ref{theorem2}}

In this section we suppose that $k=n-t-1$, where $t$ is an integer satisfying $n\ge 2t\ge 4$.

\begin{lemma}\label{lemma2-1}
Let $S,U$ be adjacent vertices of $\Gamma''_{k}(\Pi)$. Let also $M$
be the maximal singular subspace spanned by $S$ and $U$. If $T\in
\langle M]_{k}$ is a vertex of  $\Gamma''_{k}(\Pi)$ non-adjacent to
at least one of the vertices $S,U$, then there exists a vertex  $Q$
of $\Gamma''_{k}(\Pi)$ adjacent to $S,U,T$.
\end{lemma}

\proof
We suppose that $T$ is non-adjacent to $S$ (the case when $T$ is non-adjacent to $U$ is similar).
First of all, we show that the general case can be reduced to the case when
$N:=S\cap U$ and $T$ are disjoint.
Indeed, if the projective dimension of $W:=N\cap T$ is equal to $w\ge 0$
then $S,U,T$ can be naturally identified with vertices of
the graph $\Gamma''_{k-w-1}(\Pi_{W})$ which is isomorphic to
the subgraph of $\Gamma''_{k}(\Pi)$ induced on $[W\rangle_{k}$.

So, we suppose that $N$ and $T$ are disjoint.
Since $n\ge 4$, we identify $M$ with
the projective space associated to an $n$-dimensional vector space $V$.
Every non-empty singular subspace of $M$ will be identified with the corresponding subspace of
the vector space $V$.
We set
$$m:=\dim_{l} N=n-2t$$
and suppose that
$$\dim_{l}(S\cap T)=i,\;\;\;\dim_{l}(U\cap T)=j.$$
Note that $i>m$, since $S$ and $T$ are non-adjacent vertices of $\Gamma''_{k}(\Pi)$.

First we consider the case when $i+j=n-t$, i.e. $T$ is spanned by $S\cap T$ and $U\cap T$.
Since $\dim_{l}(S/N)=\dim_{l}(U/N)=t$, we can choose vectors
$$
x_{1},\dots,x_{t}\in S\setminus N\;\mbox{ and }\; y_{1},\dots,y_{t}\in U\setminus N$$
such that
$$S=N+\langle x_{1},\dots,x_{t}\rangle,\;\;\;
U=N+\langle y_{1},\dots,y_{t}\rangle,\;\;\; T=\langle x_{1},\dots,
x_{i}\rangle+\langle y_{t-j+1},\dots,y_{t}\rangle.
$$ The vectors
$x_{1}+y_{1},\dots,x_{t}+y_{t}$ are linearly independent. We define
$$Q:=N+\langle x_{1}+y_{1},\dots,x_{t}+y_{t}\rangle.$$
An easy verification shows that
$$S+Q=U+Q=T+Q=M$$
which implies that $Q$ is a vertex of $\Gamma''_{k}(\Pi)$ adjacent to $S,U,T$.

Now suppose that $$l:=n-t-(i+j)>0.$$
Then $t=i+j-m+l$. Since $i>m$, we have $i+j-m>0$.
We choose linearly independent vectors
$$x_{1},\dots,x_{i+j-m},x'_{1},\dots,x'_{l}\in S\setminus N$$
and linearly independent vectors
$$y_{1},\dots,y_{i+j-m},y'_{1},\dots,y'_{l}\in U\setminus N$$
such that
$$T=\langle x_{1},\dots,x_{i}\rangle+\langle y_{i-m+1},\dots, y_{i+j-m} \rangle+
\langle x'_{1}+y'_{1},\dots,x'_{l}+y'_{l}\rangle.$$
 We denote by $Q$ the subspace spanned by $N$ and the  vectors
$$x_{1}+y'_{1},x'_{1}+y'_{2},\dots,x'_{l-1}+y'_{l}, x'_{l}+y_{1}, x_{2}+y_{2},\dots, x_{i+j-m}+y_{i+j-m}.$$
Using the equalities
$$S=N+\langle x_{1},\dots,x_{i+j-m},x'_{1},\dots,x'_{l}\rangle,$$
$$U=N+\langle y_{1},\dots,y_{i+j-m},y'_{1},\dots,y'_{l}\rangle,$$
we establish that $S+Q=U+Q=M$.
To complete the proof we need to check that $T+Q=M$.

The conditions
$$x_{1},\dots,x_{i},y_{i-m+1},\dots, y_{i+j-m}\in T\;
\mbox{ and }\;x_{2}+y_{2},\dots, x_{i+j-m}+y_{i+j-m}\in Q$$ imply
that $x_{p}\in T+Q$ for every $p$ and $y_{p}\in T+Q$ if $p\ge 2$.
Since $x_{1}\in T$ and $x_{1}+y'_{1}\in Q$, we have $y'_{1}\in T+Q$.
Then $x'_{1}+y'_{1}\in T$ implies that $x'_{1}\in T+Q$. Step by
step, we  establish that all $x'_{q}$ and $y'_{q}$ belong to $T+Q$.
The conditions $x'_{l}\in T+Q$ and $x'_{l}+y_{1}\in Q$ guarantee
that $y_{1}\in T+Q$. Therefore, $Q+T$ coincides with $S+U=M$. $\qed$

Let $f$ be an automorphism of $\Gamma''_{k}(\Pi)$.
Show that for every $M\in {\mathcal G}_{n-1}(\Pi)$ there exists $M'\in {\mathcal G}_{n-1}(\Pi)$
such that
$$f(\langle M]_{k})\subset \langle M']_{k}.$$
We choose $S,U\in \langle M]_{k}$ such that $M$ is spanned by $S\cup
U$. Then $S$ and $U$ are adjacent vertices of $\Gamma''_{k}(\Pi)$
and the same holds for $f(S)$ and $f(U)$. Hence
$$M':=\langle f(S)\cup f(U)\rangle$$
is a maximal singular subspace.
We assert that $f(T)\in \langle M']_{k}$ for every $T\in \langle M]_{k}$
and prove this statement in several steps.

(i). If $T\in \langle M]_{k}$ is a vertex of $\Gamma''_{k}(\Pi)$
adjacent to $S$ and $U$ then $f(S)$, $f(U)$, $f(T)$ form a clique in
$\Gamma''_{k}(\Pi)$ and, by Proposition \ref{prop-cliq2}, we have
$f(T)\in \langle M']_{k}$.

(ii). Now, let $T\in \langle M]_{k}$ be a vertex of  $\Gamma''_{k}(\Pi)$
adjacent to precisely one of the vertices  $U$ and $S$.
Suppose that $T$ is adjacent to $U$ and non-adjacent to $S$.
Let $Q\in \langle M]_{k}$ be a vertex of $\Gamma''_{k}(\Pi)$ adjacent to $S,U,T$
(Lemma \ref{lemma2-1}).
By (i), $f(Q)$ belongs to $\langle M']_{k}$.
Then
\begin{equation}\label{eq2}
\langle f(U)\cup f(Q)\rangle=M'.
\end{equation}
We apply (i) to $U,Q,T$ and establish that $f(T)$ belongs to $\langle M']_{k}$.

(iii). Consider the case when $T$ is a vertex of $\Gamma''_{k}(\Pi)$
non-adjacent to both $S$ and $U$.
As above, we consider $Q\in \langle M]_{k}$ which is a vertex of $\Gamma''_{k}(\Pi)$
adjacent to $S,U,T$.
Then $f(Q)\in \langle M']_{k}$ and \eqref{eq2} holds.
We apply (ii) to $U,Q,T$ and obtain that $f(T)\in \langle M']_{k}$.

We apply the above arguments to $f^{-1}$ and
establish the existence of a bijective transformation $g$ of ${\mathcal G}_{n-1}(\Pi)$
such that
$$f(\langle M]_{k})=\langle g(M)]_{k}$$
for every $M\in {\mathcal G}_{n-1}(\Pi)$.

Let $U$ be a singular subspace of $\Pi$ such that $k<\dim_{p} U<n-1$.
Then there exist
$M_{1},M_{2}\in {\mathcal G}_{n-1}(\Pi)$ satisfying $U=M_{1}\cap M_{2}$.
We have
$$\langle U]_{k}=\langle M_{1}]_{k}\cap \langle M_{2}]_{k},$$
$$f(\langle U]_{k})=\langle g(M_{1})]_{k}\cap \langle g(M_{2})]_{k}=
\langle g(M_{1})\cap g(M_{2})]_{k}.
$$ We set
$$g(U):=g(M_{1})\cap g(M_{2})$$
and get an extension of $g$ to a transformation of
\begin{equation}\label{eq3}
{\mathcal G}_{k+1}(\Pi)\cup\dots\cup {\mathcal G}_{n-1}(\Pi)
\end{equation}
such that
$$f(\langle U]_{k})=\langle g(U)]_{k}$$
for every $U$ belonging to \eqref{eq3}.
We apply the same arguments to the pair $g^{-1},f^{-1}$
and show that $g$ is bijective.
It is clear that $f$ is inclusion preserving, i.e.
$$S\subset U\;\Longleftrightarrow\;f(S)\subset f(U)$$
for any $S,U$ belonging to \eqref{eq3}.
The latter guarantees that $f$ sends ${\mathcal G}_{i}(\Pi)$ to itself
for every integer $i$ satisfying $k<i\le n-1$.
Therefore, $f$ and $f^{-1}$ both send tops to tops
which implies that $f$ is induced by an automorphism of $\Pi$.

  \section{Valency of the relation for finite classical polar spaces}
Let $V$ be an $n$-dimensional vector space over the finite field
${\mathbb F}_{q}$ together with a non-degenerate symplectic,
hermitian or symmetric form of Witt index $d\ge 2$. For every $m\in
\{1,\dots,d\}$ we denote by ${\mathcal N}_{m}$  the set of all
$m$-dimensional totally isotropic subspaces. The associated polar
space $\Pi$ is defined as follows: the point set is ${\mathcal
N}_{1}$ and the lines are defined by elements of ${\mathcal N}_{2}$
(the line corresponding to $S\in {\mathcal N}_{2}$ consists of all
elements of ${\mathcal N}_{1}$ contained in $S$). This is a polar
space of rank $d$ and every ${\mathcal G}_{i}(\Pi)$ can be naturally
identified with ${\mathcal N}_{i+1}$.  For more information, see
\cite{taylor, wan1}.

 For $0\leq i\leq j\leq m,$ denote by $n_{i,j}$ the valency of
the relation
\begin{equation}\label{eq4}
\{(S,U)\in{\mathcal N}_{m}\times{\mathcal N}_{m}: \dim_l
(S^{\perp}\cap U)=m-i,\; \dim_l (S\cap U)=m-j\},
\end{equation}
where $S^{\perp}$ is the orthogonal complement to $S$.
If $m<d$ then \eqref{eq4} coincides with the relation ${\mathfrak R}_{i,j}$.
In the case when $m=d$,
we have $S^{\perp}\cap U=S\cap U$ for any $S,U\in {\mathcal N}_{m}$
and \eqref{eq4} is the relation ${\mathfrak R}_{i}$ defined on the dual polar space.

The valency $n_{i, j}$ is computed in \cite{Stanton} as follows:
$$
n_{i, j}=
q^{j^2+i(n-2m-2j+\frac{3}{2}i+\frac{1}{2}-\mu-\nu)}{m\brack j}{j\brack i}\\
\prod\limits_{s=0}^{j-i-1}(q^{\frac{n}{2}-m-\mu-s}-1)
(q^{\frac{n}{2}-m-\nu-s}+1)(q^{s+1}-1)^{-1},
$$
where
$\mu=\frac{1}{2}n-d$ and $\nu$ is a number such that $\mu+\nu$
equals $0,\frac{1}{2},1$ in the symplectic, hermitian and symmetric
cases, respectively.

\begin{prop}\label{prop-val}
All valencies $n_{i, j}$  are pairwise distinct.
\end{prop}

\begin{lemma}\label{5le7}
Let $f(x)$ be a polynomial of  degree at least $1$ over the rational
number field $\mathbb Q$.

{\rm (i)} If $f(x)$ has the following two factorizations:
\begin{eqnarray*}\label{4a0}
f(x)=(x^{r_1}-1)^{i_1}(x^{r_2}-1)^{i_2}\cdots
(x^{r_l}-1)^{i_l}=(x^{t_1}-1)^{j_1} (x^{t_2}-1)^{j_2}\cdots
(x^{t_h}-1)^{j_h},
\end{eqnarray*}
where $r_1>\cdots >r_l$ and $t_1>\cdots >t_h$, then $h=l$, $r_s=t_s$
and $i_s=j_s$, $s=1,\ldots,l$.

{\rm (ii)} If $f(x)$ has the following two factorizations:
\begin{eqnarray*}\label{4a0'}
f(x)=(x^r+1)(x^{r_1}-1)^{i_1}\cdots
(x^{r_l}-1)^{i_l}=(x^{t_1}-1)^{j_1}\cdots (x^{t_h}-1)^{j_h},
\end{eqnarray*}
then there   exists an $s\in \{1,\ldots,h\}$ such that $t_s=2r$.

{\rm (iii)} If $f(x)$ has the following two factorizations:
\begin{eqnarray*}\label{4a0''}
f(x)=(x^{k_1}+1)\cdots(x^{k_f}+1)(x^{r_1}-1)^{i_1}\cdots
(x^{r_l}-1)^{i_l}=(x^{t_1}-1)^{j_1}\cdots (x^{t_h}-1)^{j_h},
\end{eqnarray*}
where $k_1,\ldots,k_f$ are distinct odd  integers, then there  exist
 $s_1,\ldots s_f\in \{1,\ldots,h\}$ such that $t_{s_1}=2k_1, \ldots,
t_{s_f}=2k_f$.
\end{lemma}

\proof (i)  Consider the term  of second minimum degree of $f(x)$,
we have
$$
(-1)^{i_1+\cdots+i_l-1}i_lx^{r_l}=(-1)^{j_1+\cdots+j_h-1}j_hx^{t_h},
$$
which implies that  $j_h=i_l$ and $t_h=r_l$.  So
$$
(x^{r_1}-1)^{i_1}\cdots
(x^{r_{l-1}}-1)^{i_{l-1}}=(x^{t_1}-1)^{j_1}\cdots
(x^{t_{h-1}}-1)^{j_{h-1}},
$$
 By induction, (i) holds.

(ii) The polynomial $(x^r-1) f(x)$ has the following two
factorizations:
$$
(x^r-1) f(x)=(x^{2r}-1)(x^{r_1}-1)^{i_1}\cdots
(x^{r_l}-1)^{i_l}=(x^r-1)(x^{t_1}-1)^{j_1}\cdots (x^{t_h}-1)^{j_h}.
$$
 So (ii) holds by (i).

 The proof of (iii) is similar to that of (ii), and will be
 omitted.$\qed$

\medskip{\em Proof of Proposition \ref{prop-val}:}
Suppose $n_{a, b}= n_{a', b'}$. Since $q$ is a prime power,
we obtain
\begin{eqnarray}\label{4a1}
b^2+a(n-2m-2b+\frac{3}{2}a+\frac{1}{2}-\mu-\nu)=b'^2+a'(n-2m-2b'+\frac{3}{2}a'+\frac{1}{2}-\mu-\nu),
\end{eqnarray}
and
\begin{eqnarray}\label{4a2}
&{m\brack b}{b\brack
a}\prod\limits_{s=0}^{b-a-1}(q^{\frac{n}{2}-m-\mu-s}-1)(q^{\frac{n}{2}-m-\nu-s}+1)
(q^{s+1}-1)^{-1}\nonumber \\
=&{m\brack b'}{b'\brack
a'}\prod\limits_{s=0}^{b'-a'-1}(q^{\frac{n}{2}-m-\mu-s}-1)(q^{\frac{n}{2}-m-\nu-s}+1)
(q^{s+1}-1)^{-1}.
\end{eqnarray}

Simplifying (\ref{4a1}), we have
\begin{eqnarray}\label{4a3}
(b-b')(2b+2b'-4a')=(a'-a)(2l+3a'+3a-4b),
\end{eqnarray}
where $l=n-2m+\frac{1}{2}-\mu-\nu$. Write
$$
\begin{array} {rcl}
f_{i,j}(x)&=&\prod\limits_{s=0}^{j-i-1}(x^{n-2m-2\mu-2s}-1)(x^{n-2m-2\nu-2s}+1),\\
g_{i,j}(x)&=&\prod\limits_{s=1}^{m-j}(x^{2s}-1)\prod\limits_{s=1}^{i}(x^{2s}-1)\prod\limits_{s=1}^{j-i}(x^{2s}-1)^2.
\end{array}
$$
The equality (\ref{4a2}) implies that
$f_{a,b}(q)g_{a',b'}(q)=f_{a',b'}(q)g_{a,b}(q)$  for all prime
powers $q$, so
\begin{eqnarray}\label{4a4}
f_{a,b}(x)g_{a',b'}(x)=f_{a',b'}(x)g_{a,b}(x).
\end{eqnarray}
Computing the degree of this polynomial, we have
\begin{eqnarray}\label{4a5}
(b-b')(2l-5b-5b'+m+8a')=(a'-a)(-2l-5a-5a'+8b).
\end{eqnarray}
Equalities  (\ref{4a3}) and (\ref{4a5}) imply that
\begin{eqnarray}\label{4a6}
(b-b')(2l+2m-b-b')=(a'-a)(2l+a+a').
\end{eqnarray}

Suppose $b\neq b'.$ Without loss of generality, we assume that
$b>b'$. Since $2l+2m-b-b'>0$ and $2l+a+a'>0$, by  (\ref{4a6}) one
gets $a'>a$ and
\begin{equation}\label{ineq}
b-a-b'+a'\geq 2.
\end{equation}
Write
$$\begin{array}{rcl}
k(x)&=&\prod\limits_{s=b'-a'+1}^{b-a}(x^{2s}-1)^2,\\
h(x)&=&\prod\limits_{s=b'-a'}^{b-a-1}(x^{n-2m-2\mu-2s}-1)(x^{n-2m-2\nu-2s}+1)\prod\limits_{s=m-b+1}^{m-b'}(x^{2s}-1)
\prod\limits_{s=a+1}^{a'}(x^{2s}-1).
\end{array}
$$
The equation (\ref{4a4}) implies that $k(x)=h(x)$.

Set
$$
h_1(x)=x^{2n-4m-4\nu-4b'+4a'}-1,\quad h_2(x)=x^{2n-4m-4\nu-4b+4a+4}-1.
$$
By Lemma~\ref{5le7}  there exist
$s_1,s_2\in\{b'-a'+1,\ldots, b-a\}$ such that $h_1(x)=x^{2s_1}-1$
and $h_2(x)=x^{2s_2}-1$. It follows that
$$
s_1=n-2m-2\nu-2b'+2a',\quad s_2= n-2m-2\nu-2b+2a+2.
$$
By (\ref{ineq}) we have $s_1-s_2=2b-2a-2b'+2a'-2>b-a-(b'-a'+1),$ a
contradiction. So $b=b'$ and $a=a'$, as desired. $\qed$

\begin{rem}{\rm  It is well known that $\mathcal N_m$
has a structure of a symmetric association scheme, each relation of
which has form (\ref{eq4}). By Theorem~\ref{theorem1} every
automorphism of this scheme is induced by an automorphism of the
associated polar space $\Pi$. }\end{rem}

\section*{Acknowledgement}
W. Liu was supported by NSFC(11271004). K. Wang  was supported by
NSFC(11271047) and the Fundamental Research Funds for the Central
University of China.

\end{document}